# Varieties with maximum likelihood degree one

June Huh

ABSTRACT

We show that algebraic varieties with maximum likelihood degree one are exactly the images of reduced $A$-discriminantal varieties under monomial maps with finite fibers. The maximum likelihood estimator corresponding to such a variety is Kapranov's Horn uniformization. This extends Kapranov's characterization of $A$-discriminantal hypersurfaces to varieties of arbitrary codimension.

## 1. Main results

**1.1**

Let $X$ be a closed and irreducible subvariety of the algebraic torus

$$(\mathbb{C}^*)^m = \Big\{ \mathbf{p} = (p_1, \ldots, p_m) \in \mathbb{C}^m \mid \prod_{i=1}^m p_i \neq 0 \Big\}.$$

If $\mathbf{u} = (u_1, \ldots, u_m)$ is a set of integers, then the *likelihood function* of $X$ is defined to be

$$L = L(\mathbf{p}, \mathbf{u}) = \prod_{i=1}^m p_i^{u_i} : X \longrightarrow \mathbb{C}^*.$$

One is often interested in a statistical model $X$ contained in the hyperplane $\big\{ \sum_{i=1}^m p_i = 1 \big\}$, and in real critical points of the likelihood function corresponding to positive integer data $\mathbf{u}$. One of the critical points will provide parameters $\mathbf{p}$ which best explain the observation $\mathbf{u}$.

We refer to [CHKS06, DSS09, HKS05, PS05] for an introduction to the problem of maximum likelihood estimation in the setting of algebraic statistics. For the study of critical points of $L$ from a more geometric point of view, see [Dam99, Dam00, FK00, Huh12, OT95, Sil96, Var95].

DEFINITION 1. The *maximum likelihood degree* of $X \subseteq (\mathbb{C}^*)^m$ is the number of critical points of $L(\mathbf{p}, \mathbf{u})$ on the set of smooth points of $X$ for sufficiently general $\mathbf{u}$.

It will become clear in Section 3 that this number is finite and well-defined. We consider the following problem posed in [HKS05, Problem 14] and [Stu09, Section 3].

PROBLEM. *Find a geometric characterization of varieties with maximum likelihood degree one.*

Theorems 2 and 5 below show that the class of varieties in question is essentially the class of $A$-discriminantal varieties in the sense of Gelfand, Kapranov, and Zelevinsky [GKZ94]. In particular, there are only countably many subvarieties of $(\mathbb{C}^*)^m$ whose maximum likelihood degree is one, one for each integral matrix with $m$ columns whose column sums are zero, up to scaling of coordinates $\mathbf{p}$.

The author was partially supported by NSF grant DMS-0943832.

THEOREM 2. *A subvariety of $(\mathbb{C}^*)^m$ has maximum likelihood degree one if and only if it admits Kapranov's Horn uniformization. More precisely, the following are equivalent:*

(i) $X \subseteq (\mathbb{C}^*)^m$ *has maximum likelihood degree one.*

(ii) *There is a vector of nonzero complex constants* $\mathbf{d} = (d_1, \ldots, d_m)$, *a positive integer* $n$, *and an integral matrix*
$$B = \begin{bmatrix} b_{11} & \cdots & b_{1m} \\ \vdots & \ddots & \vdots \\ b_{n1} & \cdots & b_{nm} \end{bmatrix}$$
*whose column sums are zero, such that the rational map*
$$\Psi : \mathbb{P}^{m-1} \dashrightarrow (\mathbb{C}^*)^m, \qquad (u_1, \cdots, u_m) \longmapsto (\Psi_1, \ldots, \Psi_m),$$
*maps* $\mathbb{P}^{m-1}$ *dominantly to* $X$, *where*
$$\Psi_k(u_1, \ldots, u_m) = d_k \prod_{i=1}^{n} \Big(\sum_{j=1}^{m} b_{ij} u_j\Big)^{b_{ik}}, \qquad 1 \leqslant k \leqslant m.$$
*Here we agree that zero to the power of zero is one.*

In this case, $X \subseteq (\mathbb{C}^*)^m$ uniquely determines, and is determined by, $\Psi$.

The rational functions $\Psi_k$ are homogeneous of degree zero in the variables $\mathbf{u}$, because column sums of $B$ are assumed to be zero. The rational map $\Psi$ is the *likelihood estimator* of $X$ which maps the data vector $\mathbf{u}$ to the unique critical point of the corresponding likelihood function $L(\mathbf{p}, \mathbf{u})$.

The proof of Theorem 2 closely follows Kapranov's presentation of Horn's ideas from 1889 [Hor89]. As Kapranov remarks in [Kap91], the present paper could have been written a hundred years ago.

**1.2**

Theorem 2 shows that the set of all varieties with maximum likelihood degree one is partially ordered by taking images under monomial maps with finite fibers. To be more precise, let $B$ be an $n \times m$ integral matrix as in Theorem 2, and let $C$ be an $m \times l$ integral matrix with linearly independent rows. Consider the homomorphism
$$\phi^C : (\mathbb{C}^*)^m \longrightarrow (\mathbb{C}^*)^l, \qquad \mathbf{p} = (p_1, \ldots, p_m) \longmapsto \mathbf{p}^C := \Big(\prod_{i=1}^{m} p_i^{c_{i1}}, \ldots, \prod_{i=1}^{m} p_i^{c_{il}}\Big),$$
and the linear projection
$$\phi_C : \mathbb{P}^{l-1} \dashrightarrow \mathbb{P}^{m-1}, \qquad \mathbf{v} = (v_1, \ldots, v_l) \longmapsto C\mathbf{v} := \Big(\sum_{j=1}^{l} c_{1j} v_j, \ldots, \sum_{j=1}^{l} c_{mj} v_j\Big).$$

In the same notation, the Horn uniformization $\Psi$ of Theorem 2 can be written
$$\Psi(\mathbf{u}) = \mathbf{d} * (B\mathbf{u})^B,$$
where $*$ is the Hadamard product, the entrywise multiplication.

LEMMA 3. *For* $\mathbf{v} \in \mathbb{C}^l$ *and* $\mathbf{r}, \mathbf{d} \in (\mathbb{C}^*)^m$, *we have*
$$B(C\mathbf{v}) = (BC)\mathbf{v}, \qquad (\mathbf{r}^B)^C = \mathbf{r}^{(BC)}, \qquad (\mathbf{d} * \mathbf{r})^C = \mathbf{d}^C * \mathbf{r}^C.$$

The same rules continue to hold if **v**, **r**, **d** are replaced by matrices of appropriate sizes.

It follows that there is a commutative diagram

$$\begin{array}{ccc} \mathbb{P}^{m-1} & \dashrightarrow^{\Psi} & (\mathbb{C}^*)^m \\ \phi_C \uparrow & & \downarrow \phi^C \\ \mathbb{P}^{l-1} & \dashrightarrow_{\Psi'} & (\mathbb{C}^*)^l, \end{array}$$

where $\Psi'$ is the Horn uniformization associated to $\mathbf{d}^C$ and $BC$:

$$\Psi'(\mathbf{v}) = \mathbf{d}^C * (BC\mathbf{v})^{BC}.$$

Since $\phi_C$ is dominant and $\phi^C$ is proper, we have

$$\phi^C\left(\overline{\operatorname{im}(\Psi)}\right) = \overline{\operatorname{im}(\Psi')}.$$

COROLLARY 4. *If $X \subseteq (\mathbb{C}^*)^m$ is a closed subvariety with maximum likelihood degree one, then $\phi^C(X) \subseteq (\mathbb{C}^*)^l$ is a closed subvariety with maximum likelihood degree one.*

Note that it is necessary to assume that $C$ has rank $m$ in order to ensure that $\phi^C(X) \subseteq (\mathbb{C}^*)^l$ is closed and has maximum likelihood degree one. Note also that the maximum likelihood degrees of $X$ and $\phi^C(X)$ are different in general, even if $C$ has rank $m$. See Example 9.

**1.3**

Define a partial order on the set of all varieties with maximum likelihood degree one by

$$\left(X \subseteq (\mathbb{C}^*)^m\right) \succeq \left(X' \subseteq (\mathbb{C}^*)^l\right) \longleftrightarrow$$
$$\left(\text{there is an } m \times l \text{ integral matrix } C \text{ of rank } m \text{ such that } \phi^C(X) = X'\right).$$

The maximal elements of this partially ordered set are precisely the reduced $A$-discriminantal varieties of [GKZ94, Chapter 9], up to scaling of coordintates.

THEOREM 5. *The following are equivalent:*

(i) *$X \subseteq (\mathbb{C}^*)^m$ has maximum likelihood degree one.*
(ii) *There is a vector of nonzero complex constants $\mathbf{d} = (d_1, \ldots, d_m)$, positive integers $n$ and $k$, an integral matrix*

$$A = \begin{bmatrix} 1 & \cdots & 1 \\ a_{21} & \cdots & a_{2n} \\ \vdots & \ddots & \vdots \\ a_{k1} & \cdots & a_{kn} \end{bmatrix}$$

*whose columns generate $\mathbb{Z}^k$, and an integral matrix of rank $n - k$*

$$B = \begin{bmatrix} b_{11} & \cdots & b_{1m} \\ \vdots & \ddots & \vdots \\ b_{n1} & \cdots & b_{nm} \end{bmatrix}$$

*with $AB = 0$, such that the monomial map*

$$\mathbf{d} * \phi^B : (\mathbb{C}^*)^n \longrightarrow (\mathbb{C}^*)^m, \quad \mathbf{q} = (q_1, \ldots, q_n) \longmapsto \mathbf{d} * \mathbf{q}^B := \left(d_1 \prod_{i=1}^n q_i^{b_{i1}}, \ldots, d_m \prod_{i=1}^n q_i^{b_{im}}\right)$$

maps the $A$-discriminantal variety $\nabla_A \cap (\mathbb{C}^*)^n$ dominantly to $X$.
In this case, $\mathbf{d} * \phi^B$ factors through a monomial map with finite fibers
$$\mathbb{T}\big(\ker(A)\big) \longrightarrow (\mathbb{C}^*)^m,$$
which maps the reduced $A$-discriminantal variety in $\mathbb{T}\big(\ker(A)\big)$ birationally onto $X$.

Here $\mathbb{T}\big(\ker(A)\big) := \mathrm{Hom}\big(\ker(A), \mathbb{C}^*\big)$ is the algebraic torus whose character lattice is $\ker(A)$. If columns of $B$ form a basis of $\ker(A)$, then $X$ is the reduced $A$-discriminantal variety, up to scaling of coordinates by $\mathbf{d}$.

Our basic reference on $A$-discriminants will be [GKZ94]. The definition and basic properties of $A$-discriminantal variety and reduced $A$-discriminantal variety will be recalled in Section 3.6.

## 2. Examples and remarks

*Example* 6. A point $\{\mathbf{p}\} \in (\mathbb{C}^*)^m$ has maximum likelihood degree one. The corresponding Horn uniformization is the constant map
$$\Psi : \mathbb{P}^{m-1} \longrightarrow (\mathbb{C}^*)^m, \qquad \mathbf{u} \longmapsto \mathbf{d} * (B\mathbf{u})^B,$$
where
$$\mathbf{d} = -\mathbf{p} \qquad \text{and} \qquad B = \begin{bmatrix} 1 & \cdots & 1 \\ -1 & \cdots & -1 \end{bmatrix}.$$
The choice of $\mathbf{d}$ and $B$ is in general not unique. For example, without changing $\Psi$ one may take
$$\mathbf{d} = -\frac{27}{4}\mathbf{p} \qquad \text{and} \qquad B = \begin{bmatrix} 1 & \cdots & 1 \\ 2 & \cdots & 2 \\ -3 & \cdots & -3 \end{bmatrix}.$$
Consequently, the choice of $A$ in Theorem 5 is not unique.

*Example* 7. Consider two binary random variables, and write $\mathbf{p} = (p_{00}, p_{01}, p_{10}, p_{11})$ for the joint probabilities corresponding to four possible outcomes. The case when the two events are independent can be modeled by the algebraic variety
$$X = \Big\{ \mathbf{p} \mid p_{00}p_{11} - p_{01}p_{10} = 0, \ p_{00} + p_{01} + p_{10} + p_{11} = 1 \Big\} \subseteq (\mathbb{C}^*)^4.$$
$X$ has maximum likelihood degree one, and the likelihood function of $X$ corresponding to a given data vector $\mathbf{u} = (u_{00}, u_{01}, u_{10}, u_{11})$ is maximized at its unique critical point
$$\mathbf{p} = \Psi(\mathbf{u}) = \left( \frac{u_{0+}u_{+0}}{u_{++}^2}, \frac{u_{0+}u_{+1}}{u_{++}^2}, \frac{u_{1+}u_{+0}}{u_{++}^2}, \frac{u_{1+}u_{+1}}{u_{++}^2} \right),$$
where
$$\begin{bmatrix} u_{0+} \\ u_{1+} \\ u_{++} \\ u_{+0} \\ u_{+1} \end{bmatrix} := \begin{bmatrix} u_{00} + u_{01} \\ u_{10} + u_{11} \\ u_{00} + u_{01} + u_{10} + u_{11} \\ u_{00} + u_{10} \\ u_{01} + u_{11} \end{bmatrix}.$$
The critical point $\Psi(\mathbf{u})$ provides parameters $\mathbf{p}$ which best explains $\mathbf{u}$. Note that $\Psi$ is the Horn uniformization
$$\Psi : \mathbb{P}^3 \dashrightarrow (\mathbb{C}^*)^4, \qquad \mathbf{u} \longmapsto \mathbf{d} * (B\mathbf{u})^B,$$

with

$$\mathbf{d} = (4,4,4,4) \quad \text{and} \quad B = \begin{bmatrix} 1 & 1 & 0 & 0 \\ 0 & 0 & 1 & 1 \\ -2 & -2 & -2 & -2 \\ 1 & 0 & 1 & 0 \\ 0 & 1 & 0 & 1 \end{bmatrix}.$$

We check that $X$ is the image of an $A$-discriminant under a monomial map. Choose an integral matrix $A$ with $AB = 0$ as in Theorem 5. For example, one may take

$$A = \begin{bmatrix} 1 & 1 & 1 & 1 & 1 \\ 0 & 0 & 1 & 2 & 2 \end{bmatrix}.$$

We index the columns of $A$ by the variables $\mathbf{q} = (q_{0+}, q_{1+}, q_{++}, q_{+0}, q_{+1})$, and consider the space of all polynomials in $t$ of the form

$$F(t) = (q_{0+} + q_{1+}) + q_{++} \cdot t + (q_{+0} + q_{+1}) \cdot t^2.$$

By definition, the $A$-discriminant is the closure of the set of all such $F$ with a double root

$$\nabla_A = \left\{ \mathbf{q} \in \mathbb{C}^5 \mid q_{++}^2 - 4(q_{0+} + q_{1+})(q_{+0} + q_{+1}) = 0 \right\} \subseteq \mathbb{C}^5.$$

The monomial map

$$\mathbf{d} * \phi^B : (\mathbb{C}^*)^5 \longrightarrow (\mathbb{C}^*)^4, \quad \mathbf{q} \longmapsto \mathbf{d} * \mathbf{q}^B = \left( \frac{4q_{0+}q_{+0}}{q_{++}^2}, \frac{4q_{0+}q_{+1}}{q_{++}^2}, \frac{4q_{1+}q_{+0}}{q_{++}^2}, \frac{4q_{1+}q_{+1}}{q_{++}^2} \right)$$

maps the $A$-discriminant $\nabla_A \cap (\mathbb{C}^*)^5$ dominantly to $X$.

*Example* 8. Decomposable graphical models form an interesting class of varieties with maximum likelihood degree one. As an example, we consider the decomposable graphical model given by the following chordal graph on binary random variables $X_1, X_2, X_3, X_4$:

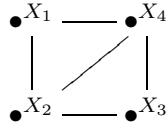

The corresponding conditional independence statement is

$$X_1 \perp\!\!\!\perp X_3 \mid \{X_2, X_4\},$$

and this defines a subvariety $X \subseteq (\mathbb{C}^*)^{16}$ cut out by four quadratic binomial equations

$$p_{1010}p_{0000} - p_{1000}p_{0010} = 0,$$
$$p_{1011}p_{0001} - p_{1001}p_{0011} = 0,$$
$$p_{1110}p_{0100} - p_{1100}p_{0110} = 0,$$
$$p_{1111}p_{0101} - p_{1101}p_{0111} = 0,$$

and one linear equation

$$\sum_{i,j,k,l} p_{ijkl} = 0.$$

$X$ has maximum likelihood degree one, and the likelihood function of $X$ corresponding to a given

data vector $\mathbf{u} = (u_{ijkl})$ is maximized at its unique critical point

$$\mathbf{p} = \Psi(\mathbf{u}) = \left( \Psi_{ijkl}(\mathbf{u}) = \frac{u_{ij+l} \cdot u_{+jkl}}{u_{++++} \cdot u_{+j+l}}, \quad 0 \leqslant i,j,k,l \leqslant 1 \right).$$

An explicit rational expression of the maximum likelihood estimator $\Psi$ for an arbitrary decomposable graphical model can be found in Lauritzen's book [Lau96, Chapter 4.4]. We invite the reader to check that the expression there indeed is a Horn uniformization.

*Example* 9. In general, the maximum likelihood degree of a variety is different from that of its image under a finite monomial map. For example, the curve

$$\{p_1^2 + p_2^2 = 1\} \subseteq (\mathbb{C}^*)^2$$

has maximum likelihood degree 4, but its image under the monomial map

$$(\mathbb{C}^*)^2 \longrightarrow (\mathbb{C}^*)^2, \qquad (p_1, p_2) \longmapsto (p_1^2, p_2^2)$$

has maximum likelihood degree 1.

*Remark* 10. Suppose $X \subseteq (\mathbb{C}^*)^m$ has maximum likelihood degree one. Then the tropicalization of $X$ can be computed from the Bergman fan of the matroid defined by the matrix $B$ of Theorem 2. See [DFS07, Section 3].

*Remark* 11. Let $\mathbb{P}^{m-1}$ be the projective space with homogeneous coordinates $p_1, \ldots, p_m$. In [HKS05, Stu09], the maximum likelihood degree is defined for a closed subvariety $X$ of

$$H := \Big\{ \mathbf{p} = (p_1, \ldots, p_m) \in \mathbb{P}^{m-1} \mid p_1 \cdots p_m (p_1 + \cdots + p_m) \neq 0 \Big\}.$$

If $\mathbf{u} = (u_1, \ldots, u_m)$ is a given data vector, then the corresponding likelihood function of $X$ is defined to be

$$L(\mathbf{p}, \mathbf{u}) = \frac{p_1^{u_1} \cdot \cdots \cdot p_m^{u_m}}{(p_1 + \cdots + p_m)^{u_1 + \cdots + u_m}} : X \longrightarrow \mathbb{C}^*.$$

We note that this setting is compatible with ours. Indeed, $H$ can be viewed as the hyperplane $\left\{ \sum_{i=1}^m p_i = 1 \right\} \subseteq (\mathbb{C}^*)^m$ by the closed embedding

$$\iota : H \longrightarrow (\mathbb{C}^*)^m, \qquad \mathbf{p} \longmapsto \left( \frac{p_1}{p_1 + \cdots + p_m}, \ldots, \frac{p_m}{p_1 + \cdots + p_m} \right).$$

The two definitions of the likelihood function of $X$ agrees under the pullback by $\iota$.

*Remark* 12. If $X \subseteq (\mathbb{C}^*)^m$ is smooth of dimension $d$, then the maximum likelihood degree of $X$ is the signed Euler-Poincaré characteristic $(-1)^d \chi(X)$. See [FK00, Huh12].

Gabber and Loeser shows in [GL96, Théorème 8.2] that a perverse sheaf is irreducible and has Euler-Poincaré characteristic one if and only if it is hypergeometric. It would be interesting to understand the relation between this result and that of the present paper. See also [LS91, LS92].

## 3. Proofs

We closely follow [GKZ94, Huh12, Kap91]. Arguments will be reproduced as needed.

**3.1**

Let $X \subseteq (\mathbb{C}^*)^m$ be a closed and irreducible subvariety of dimension $d$. We write the closed embedding by

$$\varphi : X \longrightarrow (\mathbb{C}^*)^m, \qquad \varphi = (\varphi_1, \ldots, \varphi_m).$$

For a smooth point $x$ of $X$, let $\gamma_x$ be the derivative of $\varphi$ followed by that of the left-translation $\varphi(x)^{-1} * -$

$$\gamma_x : T_x X \longrightarrow T_{\varphi(x)}(\mathbb{C}^*)^m \longrightarrow T_{\mathbf{1}}(\mathbb{C}^*)^m.$$

In local coordinates, $\gamma_x$ is represented by the logarithmic jacobian matrix

$$\left( \frac{\partial \log \varphi_i}{\partial x_j} \right), \quad 1 \leqslant i \leqslant m, \quad 1 \leqslant j \leqslant d.$$

This defines the logarithmic Gauss map to the Grassmannian of the Lie algebra $\mathfrak{g}$ of $(\mathbb{C}^*)^m$

$$\gamma : X \dashrightarrow \mathrm{Gr}(d, \mathfrak{g}), \qquad x \longmapsto \mathrm{im}(\gamma_x).$$

*Remark* 13. When $X$ is a hypersurface, $\gamma$ agrees with the logarithmic Gauss map of [GKZ94, Section 9.3]:

$$\gamma : X \dashrightarrow \mathrm{Gr}(m-1, \mathfrak{g}) = \mathbb{P}(\mathfrak{g}^\vee).$$

$X$ has maximum likelihood degree one if and only if $\gamma$ is birational, because the set of critical points of the likelihood function of $X$ corresponding to $\mathbf{u} = (u_1, \ldots, u_m)$ is the fiber of $\gamma$ over the point

$$\sum_{i=1}^m u_i \cdot \mathrm{dlog}(p_i) \in H^0\Big( (\mathbb{C}^*)^m, \Omega^1_{(\mathbb{C}^*)^m} \Big) \simeq \mathfrak{g}^\vee.$$

Kapranov states in [Kap91, Theorem 1.3] that

- if $X$ is a reduced $A$-discriminantal hypersurface, then $\gamma$ is birational, and
- if $\gamma$ is birational, then $X$ is a reduced $A$-discriminantal hypersurface, up to an automorphism of the ambient torus.

As pointed out in [CD07, Section 2], a small correction needs to be made in the second statement. If $\gamma$ is birational, then there is a monomial map with finite fibers

$$\mathbb{T}\big(\ker(A)\big) \longrightarrow (\mathbb{C}^*)^m,$$

which maps the reduced $A$-discriminantal variety in $\mathbb{T}\big(\ker(A)\big)$ birationally onto $X$.

**3.2**

We write $\mathbf{p} = (p_1, \ldots, p_m)$ for the coordinate functions of $(\mathbb{C}^*)^m$ as before. This defines a basis of the dual $\mathfrak{g}^\vee \simeq \mathbb{C}^m$ corresponding to differential forms

$$\mathrm{dlog}(p_1), \ldots, \mathrm{dlog}(p_m) \in H^0\Big( (\mathbb{C}^*)^m, \Omega^1_{(\mathbb{C}^*)^m} \Big).$$

Hereafter we fix this choice of basis of $\mathfrak{g}^\vee$, and identify $\mathfrak{g}^\vee$ with the space of data vectors $\mathbf{u}$. Consider the vector bundle homomorphism defined by the pullback of differential forms

$$\gamma^\vee : X_{\mathrm{sm}} \times \mathfrak{g}^\vee \longrightarrow \Omega^1_{X_{\mathrm{sm}}}, \qquad (x, \mathbf{u}) \longmapsto \sum_{i=1}^m u_i \cdot \mathrm{dlog}(\varphi_i)(x), \qquad \mathbf{u} = (u_1, \ldots, u_m).$$

The induced linear map $\gamma_x^\vee$ between the fibers over a smooth point $x$ is dual to the injective linear map of the previous subsection

$$\gamma_x : T_x X \longrightarrow \mathfrak{g}.$$

Therefore $\gamma^\vee$ is surjective and $\ker(\gamma^\vee)$ is a vector bundle.

DEFINITION 14. The *variety of critical points* of $X \subseteq (\mathbb{C}^*)^m$ is defined to be the closure

$$\mathfrak{X} := \overline{\mathbb{P}\big(\ker(\gamma^\vee)\big)} \subseteq X \times \mathbb{P}(\mathfrak{g}^\vee).$$

Note that $\ker(\gamma^\vee)$ is a vector bundle of rank $m - d$. Therefore $\mathfrak{X}$ is irreducible and

$$\dim \mathfrak{X} = m - 1.$$

If $\mathbf{u}$ is integral and $x$ is a smooth point of $X$, then $(x, \mathbf{u})$ is in $\mathfrak{X}$ if and only if $x$ is a critical point of the likelihood function of $X$ corresponding to $\mathbf{u}$. Since $\dim \mathfrak{X} = \dim \mathbb{P}(\mathfrak{g}^\vee)$, the maximum likelihood degree of $X$ is finite and well-defined.

**3.3**

Write $\mathbb{P}^{m-1}$ for the projective space $\mathbb{P}(\mathfrak{g}^\vee)$ with the homogeneous coordinates $\mathbf{u}$. Let $\Psi$ be a rational map

$$\Psi : \mathbb{P}^{m-1} \dashrightarrow (\mathbb{C}^*)^m, \qquad \Psi = (\Psi_1, \ldots, \Psi_m).$$

Each component of $\Psi$ should be a homogeneous rational function of degree zero in the variables $\mathbf{u}$. We have Euler's relation

$$\sum_{i=1}^m u_i \frac{\partial \log \Psi_j}{\partial u_i} = 0, \qquad 1 \leqslant j \leqslant m.$$

The following lemma will play a central role in the proof of Theorem 2.

LEMMA 15. *Suppose that the closure of the image of $\Psi$ is $X$. Then the following conditions are equivalent.*

 (i) *$\mathfrak{X}$ is the closure of the graph of $\Psi$.*
 (ii) *The graph of $\Psi$ is contained in $\mathfrak{X}$.*
 (iii) *We have*

$$\sum_{i=1}^m u_i \frac{\partial \log \Psi_i}{\partial u_j} = 0, \qquad 1 \leqslant j \leqslant m.$$

 (iv) *We have*

$$\frac{\partial \log \Psi_i}{\partial u_j} = \frac{\partial \log \Psi_j}{\partial u_i}, \qquad 1 \leqslant i \leqslant m, \quad 1 \leqslant j \leqslant m.$$

*Proof.* Since $\mathfrak{X}$ is irreducible of dimension $m - 1$, (i) and (ii) are equivalent. We prove that (ii) and (iii) are equivalent.

By generic smoothness, for a sufficiently general $\mathbf{u} \in \mathfrak{g}^\vee$,

– $\Psi(\mathbf{u})$ is a smooth point of $X$, and
– $\Psi|_U : U \to X$ is a submersion for a small neighborhood $U$ of $\mathbf{u}$ in $\mathbb{P}^{m-1}$.

Note from the construction of $\mathfrak{X}$ that the graph of $\Psi$ is contained in $\mathfrak{X}$ if and only if such $\mathbf{u}$ is contained in the kernel of $\gamma^{\vee}_{\Psi(\mathbf{u})}$. Dually, this condition is satisfied if and only if the hyperplane of $\mathfrak{g}$ defined by $\mathbf{u}$ contains the image of

$$\gamma_{\Psi(\mathbf{u})} : T_{\Psi(\mathbf{u})} X \longrightarrow \mathfrak{g}.$$

We express this last condition in terms of equations.

Fix a sufficiently general $\mathbf{u} \in \mathfrak{g}^{\vee}$ as above. The key player is the linear mapping

$$\Phi : \mathfrak{g}^{\vee} \longrightarrow \mathfrak{g},$$

defined as the composition

$$\mathfrak{g}^{\vee} \simeq T_{\mathbf{u}}\mathfrak{g}^{\vee} \longrightarrow T_{\mathbf{u}}\mathbb{P}^{m-1} \longrightarrow T_{\Psi(\mathbf{u})}(\mathbb{C}^*)^m \longrightarrow \mathfrak{g}.$$

The first is the derivative of the quotient map defining $\mathbb{P}^{m-1}$, the second is the derivative of $\Psi$, and the last is the derivative of the left-translation $\Psi(\mathbf{u})^{-1} * -$. In coordinates, $\Phi$ is represented by the logarithmic jacobian matrix

$$\begin{bmatrix} \frac{\partial \log \Psi_1}{\partial u_1} & \cdots & \frac{\partial \log \Psi_1}{\partial u_m} \\ \vdots & \ddots & \vdots \\ \frac{\partial \log \Psi_m}{\partial u_1} & \cdots & \frac{\partial \log \Psi_m}{\partial u_m} \end{bmatrix}.$$

By the genericity assumption on $\mathbf{u}$ made above, the columns of the logarithmic jacobian matrix generate the image of $\gamma_{\Psi(\mathbf{u})}$ in $\mathfrak{g}$. Therefore the image of $\gamma_{\Psi(\mathbf{u})}$ is contained in the hyperplane defined by $\mathbf{u}$ if and only if

$$\begin{bmatrix} \frac{\partial \log \Psi_1}{\partial u_1} & \cdots & \frac{\partial \log \Psi_m}{\partial u_1} \\ \vdots & \ddots & \vdots \\ \frac{\partial \log \Psi_1}{\partial u_m} & \cdots & \frac{\partial \log \Psi_m}{\partial u_m} \end{bmatrix} \begin{bmatrix} u_1 \\ \vdots \\ u_m \end{bmatrix} = 0.$$

This proves the equivalence of (ii) and (iii).

Now suppose that (iii) holds. Then

$$\frac{\partial}{\partial u_j}\left(\sum_{k=1}^m u_k \log \Psi_k\right) = \log \Psi_j + \sum_{i=1}^m u_i \frac{\partial \log \Psi_i}{\partial u_j} = \log \Psi_j,$$

and hence

$$\frac{\partial \log \Psi_j}{\partial u_i} = \frac{\partial}{\partial u_i}\frac{\partial}{\partial u_j}\left(\sum_{k=1}^m u_k \log \Psi_k\right) = \frac{\partial}{\partial u_j}\frac{\partial}{\partial u_i}\left(\sum_{k=1}^m u_k \log \Psi_k\right) = \frac{\partial \log \Psi_i}{\partial u_j}.$$

Therefore (iii) implies (iv). Lastly, (iii) is obtained from Euler's relation and (iv). $\square$

**3.4**

We continue to assume that $\Psi$ is a rational function from $\mathbb{P}^{m-1}$ to $(\mathbb{C}^*)^m$ whose components are homogeneous of degree zero in the variables $\mathbf{u}$. The following statement can be found in [Kap91, Proposition 3.1], where Kapranov attributes the result to Horn [Hor89].

LEMMA 16. *The following conditions are equivalent.*

(i) We have

$$\frac{\partial \log \Psi_i}{\partial u_j} = \frac{\partial \log \Psi_j}{\partial u_i}, \qquad 1 \leqslant i \leqslant m, \quad 1 \leqslant j \leqslant m.$$

(ii) There is a vector of nonzero constants $\mathbf{d} = (d_1, \ldots, d_m)$, a positive integer $n$, and an integral matrix
$$B = \begin{bmatrix} b_{11} & \cdots & b_{1m} \\ \vdots & \ddots & \vdots \\ b_{n1} & \cdots & b_{nm} \end{bmatrix}$$
whose column sums are zero, such that
$$\Psi_k(u_1, \ldots, u_m) = d_k \prod_{i=1}^{n} \Big(\sum_{j=1}^{m} b_{ij} u_j\Big)^{b_{ik}}, \qquad 1 \leqslant k \leqslant m.$$
Here we agree that zero to the power of zero is one.

*Proof that (i) implies (ii).* We employ the notation introduced in Lemma 3. Use unique factorization in $\mathbb{C}[u_1, \ldots, u_m]$ to write
$$\Psi = \mathbf{f}^B,$$
where
- $\mathbf{f} = (f_1, \ldots, f_n)$ is a vector of irreducible homogeneous polynomials of degrees $(\delta_1, \ldots, \delta_n)$ in the variables $\mathbf{u}$, and
- $B$ is an $n \times m$ integral matrix such that
$$\begin{bmatrix} \delta_1, \ldots, \delta_n \end{bmatrix} \begin{bmatrix} b_{11} & \cdots & b_{1m} \\ \vdots & \ddots & \vdots \\ b_{n1} & \cdots & b_{nm} \end{bmatrix} = 0.$$

We may assume that $f_i$ and $f_j$ are relatively prime to each other for $i \neq j$. Now (i) reads
$$\sum_{k=1}^{n} \Big(b_{ki} \frac{\partial f_k}{\partial u_j} - b_{kj} \frac{\partial f_k}{\partial u_i}\Big) f_1 \cdots \hat{f}_k \cdots f_n = 0.$$
Since the polynomial inside the parenthesis has degree one less than $f_k$, which is relatively prime to all the other components of $\mathbf{f}$, we have
$$b_{ki} \frac{\partial f_k}{\partial u_j} - b_{kj} \frac{\partial f_k}{\partial u_i} = 0.$$
Therefore there are homogeneous polynomials $g_k$ in $\mathbf{u}$ such that
$$\frac{\partial f_k}{\partial u_i} = b_{ki} \cdot g_k.$$
Now use Euler's relation to note that
$$\delta_k f_k = \sum_{i=1}^{m} u_i \frac{\partial f_k}{\partial u_i} = \Big(\sum_{i=1}^{m} b_{ki} u_i\Big) g_k.$$
Since $f_k$ are assumed to be irreducible, $g_k$ should be nonzero constants. This shows that
$$\mathbf{f} = \mathbf{e} * B\mathbf{u}$$
for a vector of nonzero constants $\mathbf{e} = (e_1, \ldots, e_n)$. The proof is completed by setting
$$\mathbf{d} = \mathbf{e}^B.$$
□

### 3.5

*Proof of Theorem 2.* Suppose that $X$ has maximum likelihood degree one. Let $\mathrm{pr}_1$ and $\mathrm{pr}_2$ be the projections

$$\begin{array}{ccc} & \mathfrak{X} & \\ {}_{\mathrm{pr}_1}\swarrow & & \searrow_{\mathrm{pr}_2} \\ X & & \mathbb{P}^{m-1} \end{array}$$

The assumption made on $X$ is equivalent to the statement that $\mathrm{pr}_2$ is a birational morphism. Let $\mathrm{pr}_2^{-1}$ be the rational inverse of $\mathrm{pr}_2$, and define

$$\Psi := \mathrm{pr}_1 * \mathrm{pr}_2^{-1} : \mathbb{P}^{m-1} \dashrightarrow (\mathbb{C}^*)^m.$$

Since the graph of $\Psi$ is contained in $\mathfrak{X}$, Lemma 15 and Lemma 16 prove what we want.

Conversely, suppose that $\Psi$ is a rational map of the form

$$\Psi = \mathbf{d} * (B\mathbf{u})^B,$$

which maps dominantly to $X$. By Lemma 15 and Lemma 16, $\mathfrak{X}$ is the closure of the graph of $\Psi$. This shows that $\mathrm{pr}_2$ is a birational morphism.

The above argument also shows that $X \subseteq (\mathbb{C}^*)^m$ uniquely determines, and is determined by, the rational map $\Psi$. $\square$

### 3.6

Before proceeding to the proof of Theorem 5, we recall the definition and basic properties of $A$-discriminantal varieties and reduced $A$-discriminantal varieties, following [GKZ94, Chapter 9]. Some notations are adjusted for the internal consistency of the present paper.

Let $A$ be an integral matrix of the form

$$A = \begin{bmatrix} 1 & \cdots & 1 \\ a_{21} & \cdots & a_{2n} \\ \vdots & \ddots & \vdots \\ a_{k1} & \cdots & a_{kn} \end{bmatrix}$$

whose columns generate $\mathbb{Z}^k$. Write $\{\omega_1, \ldots, \omega_n\}$ for the set of column vectors of $A$, and consider the affine space $\mathbb{C}^n$ of Laurent polynomials of the form

$$F(\mathbf{t}) = \sum_{i=1}^n q_i \cdot \mathbf{t}^{\omega_i}, \qquad \mathbf{q} = (q_1, \ldots, q_n) \in \mathbb{C}^n, \qquad \mathbf{t} = (t_1, \ldots, t_k).$$

DEFINITION 17. The *$A$-discriminantal variety* $\nabla_A$ is the closure of the set

$$\nabla_A^* = \left\{ F \in \mathbb{C}^n \mid \{F = 0\} \text{ has a singular point in } (\mathbb{C}^*)^k \right\} \subseteq \mathbb{C}^n.$$

The projective dual of $\mathbb{P}(\nabla_A) \subseteq \mathbb{P}^{n-1}$ is the toric variety $X_A \subseteq \mathbb{P}^{n-1}$, defined as the closure of the image of the monomial map

$$(\mathbb{C}^*)^k \longrightarrow \mathbb{P}^{n-1}, \qquad \mathbf{t} \longmapsto \mathbf{t}^A = (\mathbf{t}^{\omega_1}, \ldots, \mathbf{t}^{\omega_n}).$$

Let $\mathscr{B}$ be an integral matrix whose columns form a basis of $\ker A$. In other words, $\mathscr{B}$ is a Gale dual of $A$. We have exact sequences

$$0 \longrightarrow \ker(A) \simeq \mathbb{Z}^{n-k} \xrightarrow{\mathscr{B}} \mathbb{Z}^n \xrightarrow{A} \mathbb{Z}^k \longrightarrow 0$$

and
$$0 \longrightarrow (\mathbb{C}^*)^k \xrightarrow{\phi^A} (\mathbb{C}^*)^n \xrightarrow{\phi^{\mathscr{B}}} (\mathbb{C}^*)^{n-k} \simeq \mathbb{T}\big(\ker(A)\big) \longrightarrow 0.$$

Note that $\nabla_A$ is invariant under the action of $(\mathbb{C}^*)^k$.

DEFINITION 18. The *reduced A-discriminantal variety* $\widetilde{\nabla}_A$ is the image of $\nabla_A \cap (\mathbb{C}^*)^n$ in $\mathbb{T}(\ker A)$.

Reduced $A$-discriminantal varieties admit a Horn uniformization [GKZ94, Theorem 9.3.3]:

THEOREM 19. *Let $\mathscr{P}$ be the Horn uniformization*
$$\mathscr{P} : \mathbb{P}^{n-k-1} \dashrightarrow (\mathbb{C}^*)^{n-k}, \qquad \mathbf{v} \longmapsto (\mathscr{B}\mathbf{v})^{\mathscr{B}}.$$
*Then the closure of the image of $\Psi$ is the reduced A-discriminantal variety $\widetilde{\nabla}_A$.*

### 3.7

*Proof of Theorem 5.* Suppose that $X$ has maximum likelihood degree one. Then, by Theorem 2, there is a set of nonzero constants $\mathbf{d} = (d_1, \ldots, d_m)$ and an $n \times m$ integral matrix $B$ whose column sums are zero such that the Horn uniformization
$$\Psi : \mathbb{P}^{m-1} \dashrightarrow (\mathbb{C}^*)^m, \qquad \mathbf{u} \longmapsto \mathbf{d} * (B\mathbf{u})^B$$
maps dominantly to $X$. Write $n-k$ for the rank of $B$, and consider the largest subgroup $\mathbb{Z}^n$ of rank $n-k$ containing all the columns of $B$. Let $\mathscr{B}$ be a matrix whose columns form a basis of this subgroup. Let $A$ and $C$ be integral matrices such that

- $AB = 0$,
- $B = \mathscr{B}C$,
- the first row of $A$ is $(1, \ldots, 1)$, and
- the top row of the diagram below is exact:

$$\begin{array}{ccccccccc}
0 & \longrightarrow & \mathbb{Z}^{n-k} & \xrightarrow{\mathscr{B}} & \mathbb{Z}^n & \xrightarrow{A} & \mathbb{Z}^k & \longrightarrow & 0 \\
& & C \uparrow & \nearrow B & & & & & \\
& & \mathbb{Z}^m & & & & & &
\end{array}$$

Let $\mathscr{P}$ be the Horn uniformization
$$\mathscr{P} : \mathbb{P}^{n-k-1} \dashrightarrow (\mathbb{C}^*)^{n-k}, \qquad \mathbf{v} \longmapsto (\mathscr{B}\mathbf{v})^{\mathscr{B}}.$$

In the notation introduced in Section 1.2, we have a commutative diagram

$$\begin{array}{ccccc}
\mathbb{P}^{n-k-1} & \xdashrightarrow{\mathscr{P}} & (\mathbb{C}^*)^{n-k} & \xleftarrow{\phi^{\mathscr{B}}} & (\mathbb{C}^*)^n \\
\phi_C \uparrow & & \mathbf{d}*\phi^C \downarrow & \swarrow \mathbf{d}*\phi^B & \\
\mathbb{P}^{m-1} & \xdashrightarrow{\Psi} & (\mathbb{C}^*)^m & &
\end{array}$$

By Theorem 19, the commutative diagram restricts to that of dominant mappings

$$\begin{array}{ccc}
\mathbb{P}^{n-k-1} & \dashrightarrow \widetilde{\nabla}_A & \longleftarrow \nabla_A \cap (\mathbb{C}^*)^n \\
\uparrow & \downarrow & \swarrow \\
\mathbb{P}^{m-1} & \dashrightarrow X &
\end{array}$$

This proves that (i) implies (ii). Commutativity of the above diagrams also show that the monomial map with finite fibers

$$\mathbf{d} * \phi^C : (\mathbb{C}^*)^{n-k} \longrightarrow (\mathbb{C}^*)^m$$

restricts to a birational isomorphism

$$\widetilde{\nabla}_A \longrightarrow X.$$

Indeed, by Lemma 15, a fiber of $\Psi$ over a general point of $X$ is connected.

Conversely, suppose that $X$ satisfies the condition (ii). Theorem 19 and Theorem 2 show that a reduced $A$-discriminantal variety has maximum likelihood degree one. Therefore, by Corollary 4, $X$ has maximum likelihood degree one. □

## Acknowledgements

The author thanks Mikhail Kapranov and Bernd Sturmfels for the suggestion of the problem and for helpful discussions. He is grateful to Alicia Dickenstein, Seth Sullivant, and the anonymous referee for useful comments.

June Huh   junehuh@umich.edu

Department of Mathematics, University of Michigan, Ann Arbor, MI 48109, USA